\documentclass[10pt,a4paper]{article}
\usepackage[latin1]{inputenc}
\usepackage{amsmath}
\usepackage{amsfonts}
\usepackage{amssymb}
\usepackage{graphicx}
\begin{document}
\bibliographystyle{plain}
\title{An Explicit Formula for the Spherical Curves with Constant Torsion}
\author{Demetre Kazaras \\ University of Oregon
\and Ivan Sterling \\ St. Mary's College of Maryland}
\date{}
\maketitle

\section{Introduction}
\begin{figure}
\begin{center}
\includegraphics[scale=.6]{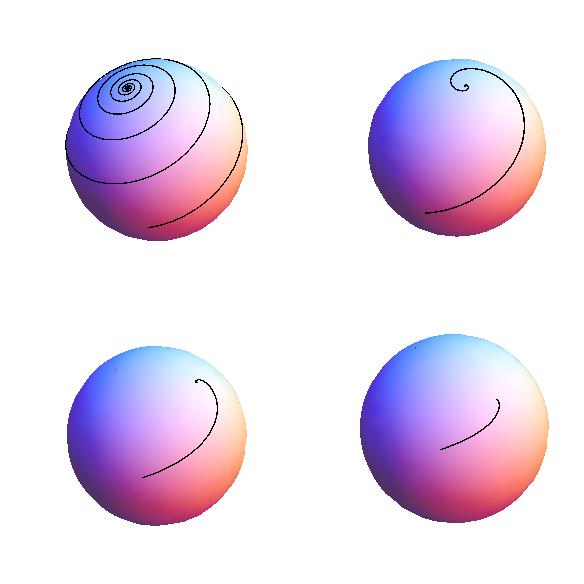}
\caption{The curves of torsion $\tau=.1,.5,1,2$ on the unit 2-sphere}\label{peter}
\end{center}
\end{figure}

The purpose of this article is to give an explicit formula for all curves of constant torsion $\tau$ in the unit two-sphere $S^2(1)$.  These curves and their basic properties have been known since the 1890's, and some of these properties are discussed in the Appendix.  Some example curves, computed with a standard ODE package, with $\tau=.1,.5,1,2$ are shown in Figure \ref{peter}.  Though their existence and some of their general properties were known, our explicit formulas for them, in terms of hypergeometric functions, are new.  

Curves of constant torsion are also of interest because all asymptotic curves on any pseudo-spherical surface (that is, surfaces in $\mathbb{R}^3$ with constant negative Gauss curvature) are of constant torsion.  Furthermore, any pair of curves with constant torsion $\pm \tau$, intersecting at one point, define an essentially unique pseudo-spherical surface.  A complete classification of curves of constant torsion in $\mathbb{R}^3$, in the context of integrable geometry, is a work in progress and is related to the corresponding unfinished classification of pseudo-spherical surfaces. 

The authors would like to thank the referees and the editor for their suggestions to improve the original version of this paper.

\section{General Setting} \label{two}
\subsection{General Curves in $\mathbb{R}^3$}
Let 
\[\gamma:(a,b) \longrightarrow \mathbb{R}^3\]
be a regular (i.e. nonzero speed) $C^\infty$ curve in $\mathbb{R}^3$ with nonzero curvature.  The speed $v$, curvature $\kappa$ and torsion $\tau$ of $\gamma$ are given by:
\[v=\Vert \gamma' \Vert,\;\;
\kappa = \frac{\Vert \gamma' \times \gamma'' \Vert}{\Vert \gamma' \Vert^3},\;\;
\tau = \frac{[\gamma' \gamma'' \gamma''']}{\Vert \gamma' \times \gamma'' \Vert^2}. \]
The unit tangent $T$ is given by
\begin{equation} \label{tangent}T=\frac{\gamma'}{v}.\;\; \end{equation}
The unit normal and unit binormal are given by
\begin{equation*} N=\frac{T'}{v \kappa},\;\; B=T \times N. \end{equation*}
These are related by the Frenet formulas:
\begin{equation} \label{frenet}
\begin{array}{ccccc} 
T' &=&&  v \kappa N  &\\
N' &=& - v \kappa T &&+ v \tau B   \\
B' &=&&-v \tau N & 
\end{array}
\end{equation}
Curves $\gamma$ with prescribed differentiable curvature $\kappa >0$ and torsion $\tau$ can be found by integrating (\ref{frenet}) and (\ref{tangent}).  Up to re-parametrization (see below) a curve in $\mathbb{R}^3$ is determined, up to a rigid motion of $\mathbb{R}^3$, by its curvature $\kappa$ and torsion $\tau$.

\subsection{Changing Parametrizations}
Given $\gamma(t)$, then the arc-length function $s(t)$ of $\gamma(t)$ is given by 
\[s(t)=\int_a^t \Vert \gamma'(u)\Vert \; du. \]
Note that since $s(t)$ is increasing, it has an inverse $t(s)$.  To obtain a unit-speed reparametrization $\gamma^{unit}$ of $\gamma$ we let 
\[\gamma^{unit}(s) := \gamma(t(s))\]
We denote the parameter of a unit-speed curve by the letter $s$. 

On the other hand, if we are given a unit-speed curve $\gamma^{unit}(s)$, we may wish to find a reparametrization $\gamma$ of $\gamma^{unit}$ by letting $t(s)$ be some special monotone function.   In that case we have
\[\gamma(t):=\gamma^{unit}(s(t))\]
where $s(t)$ is the inverse of $t(s)$.
\subsection{Spherical Curves}
$\gamma$ is called a spherical curve if (for some $r>0$) $\gamma(t) \in S^2(r) \; \forall t$.  The speed, curvature and torsion of a spherical curve satisfy \cite{G}

\begin{equation} \label{sphereeqn}
\kappa^2 \tau^2(\kappa^2 r^2 -1)=\kappa'^2 v^2.
\end{equation}

\subsection{Effect of Homothety on Curvature and Torsion}
If $\tilde{\gamma}(t):=\lambda \gamma(t)$, then $\tilde{\kappa}(t)=\frac{\kappa(t)}{\lambda}$ and $\tilde{\tau}(t)=\frac{\tau(t)}{\lambda}$.  Thus, a curve of constant torsion $\tau_1$ on a sphere of radius $r_1$ corresponds by homothety to a curve of constant torsion $\tau_2=\frac{r_1}{r_2} \tau_1$ on a sphere of radius $r_2$.  

In other words, any spherical curve of constant positive torsion corresponds to precisely one spherical curve with $\tau=1$ as well as to precisely one curve of constant positive torsion on the unit sphere.  Without loss of generality, we consider only spherical curves of constant positive torsion on the unit sphere.

\subsection{Constant Torsion Unit-Speed Curves on the Unit Sphere}
If $r=1$,  $\tau$ is a positive constant, and  $\gamma: (a,b) \rightarrow S^2(1)$ is unit-speed, then equation (\ref{sphereeqn}) is an ordinary differential equation in $\kappa$ (notice  $\kappa \geq 1$ holds for any curve on the unit sphere):
\begin{equation} \label{grayode} \kappa'^2=\kappa^2 \tau^2 (\kappa^2-1). \end{equation}
The  general solution to equation (\ref{grayode}) is given by:
\begin{equation} \label{spherekappa} \kappa = \csc(\tau s + C),
\;\; \frac{-C}{\tau} < s < \frac{-C+\pi}{\tau}. \end{equation}
Notice we use the parameter $s$ instead of $t$ since $\gamma$ is a unit-speed curve.  Furthermore, $\kappa(s)$ is decreasing on $(\frac{-C}{\tau}, \frac{-C+\frac{\pi}{2}}{\tau})$.

\subsection{Our Goal}
As mentioned,  a unit-speed curve $\gamma$ is determined up to rigid motion by its curvature and torsion.  However, in general it is not possible to explicitly solve for $\gamma$ given $\kappa>0$ and $\tau > 0$.  Spherical curves of constant torsion provide an interesting and natural example to study.  They were considered by classical geometers and the formula (\ref{spherekappa}) was known.  Even though the formula for $\kappa$ is so simple, no explicit solutions for $\gamma$ were found.  This is most likely because the integration methods that we found neccessary were not developed until decades later.  By choosing a special re-parametrization and using functions defined in the 1940's, we were successful in obtaining an explicit formula for $\gamma$ involving hypergeometric functions.

\section{Explicit Formulas}
\subsection{The Radius of Curvature Parametrization}
In curve theory parametrization by the curvature is called the ``natural parametrization".  In our case, when the natural parametrization is used, the domain of definition lies outside the radius of convergence of the resulting hypergeometric solutions.   To avoid having to deal with the problem of analytically continuing hypergeometric functions beyond their radii of convergence we instead parametrize by the reciprocal of the curvature, which is called the radius of curvature.
 
We seek unit-speed curves $\gamma^{unit}:(\frac{-C}{\tau},\frac{-C+\frac{\pi}{2}}{\tau}) \rightarrow S^2(1)$ of constant torsion $\tau >0$ on the unit sphere.  In order to simplify the Frenet equations, we reparametrize $\gamma^{unit}$ by $t(s)=\frac{1}{\kappa(s)}=\sin(\tau s + C)$. Since $\frac{1}{\kappa(s)}$ is increasing on its domain, the inverse $s(t)$, $s:(0,1) \longrightarrow (\frac{-C}{\tau},\frac{-C+\frac{\pi}{2}}{\tau})$, exists and we have
\[\gamma(t)=\gamma^{unit}(s(t))=\gamma^{unit}(\frac{\sin^{-1}(t)-C}{\tau}),
\;\; 0<t<1 .\]
One can recover $\gamma^{unit}$ from $\gamma$ by reversing the process.  Note that
\[v= \Vert \gamma' \Vert = \Vert \gamma^{unit'} \Vert |s'(t)|=|s'(t)|
=\frac{1}{\tau  \sqrt{1-t^2}}.\]
With $\kappa=\frac 1t$ the Frenet equations (\ref{frenet}) become

\begin{subequations}
\label{natfrenet}
\begin{alignat}{3}
T' &= & \frac{v  N}{t}  &\label{natfrenet1}  \\ 
N' &=\frac{- v T}{t} &&+ v \tau B  \label{natfrenet2} \\ 
B' &= & -v \tau N &  \label{natfrenet3}
\end{alignat}
\end{subequations}

Recall $\gamma' = v T$.  Thus as a preliminary step we will compute $T$.  Namely, we want to solve (\ref{natfrenet}) for $T$.  

From (\ref{natfrenet1}) and (\ref{natfrenet2}) we have 
\begin{equation} N=t \sqrt{1-t^2} \tau T', \label{Neqn} \end{equation}
and
\begin{equation} B=\sqrt{1-t^2} N'+\frac{1}{t \tau} T. \label{Beqn} \end{equation}
(\ref{Neqn}) and (\ref{natfrenet3}) yield
\begin{equation} B'=-\tau t T'. \label{Bprime1} \end{equation}
On the other hand differentiating (\ref{Neqn}) yields
\begin{equation} N'= \frac{\tau}{\sqrt{1-t^2}}(t (1-t^2) T''+ (1-2 t^2) T'). \label{Nprime} \end{equation}
Plugging this into (\ref{Beqn}) yields
\begin{equation} B=\frac{-1}{t \tau}(t^2(t^2-1) \tau^2 T''+t(2 t^2-1) \tau^2 T'-T). \label{B} \end{equation}
Hence
\begin{equation} B'=\frac{-1}{t^2 \tau}(t^3(t^2-1)\tau^2T'''+t^2(5t^2-2)\tau^2T''+t(4t^2 \tau ^2-1)T'+T). \label{Bprime2} \end{equation}
Equating (\ref{Bprime1}) and (\ref{Bprime2}) and simplifying we arrive at
\begin{equation} \label{ODET}
t^3(t^2-1)\tau^2T'''+t^2(5t^2-2)\tau^2T''+t(3t^2 \tau ^2-1)T'+T=0.
\end{equation}
This is a third order linear homogeneous differential equations with non-constant coefficients.  In general it is not possible to find a closed form solution for such an equation.  However, this is one of the special cases where one can find hypergeometric type solutions.  These methods were developed in the 1940's, and hence were not available to the classical (1890's) geometers.

\subsection{Initial Conditions}
To arrive at initial conditions for our ODE, we find initial conditions for $T$,  $N$, and $B$ and use the Frenet equations (\ref{natfrenet}) to arrive at initial conditions for $T$, $T'$, and $T''$.  We let $T=(T_1,T_2,T_3)$, $N=(N_1,N_2,N_3)$, and $B=(B_1,B_2,B_3)$.  For $V=T,N,B,T', \mbox{or}\; T''$ we us the notation $V_{i_0} := V_i(t_0)$.

$T$ and $N$ are unit vectors ($\Vert T \Vert =1$ and $\Vert N\Vert=1$) so we have
\[ |T_{1_0}| \leq 1,\;\;|T_{2_0}|\leq \sqrt{1-T_{1_0}},\;\;T_{3_0}=\sqrt{1-T_{1_0}^2-T_{2_0}^2}\]
\[ |N_{1_0}|\leq 1,\;\;|N_{2_0}|\leq \sqrt{1-N_{1_0}},\;\;N_{3_0}=\sqrt{1-N_{1_0}^2-N_{2_0}^2}\]
Also $T$ is  orthogonal to $N$ ($T \cdot N = 0$) 
\[N_{1_0}T_{1_0}+N_{2_0}T_{2_0}+N_{3_0}T_{3_0}=0.\]
By $B=T\times N,$ we have
\[B_{1_0}=T_{2_0}N_{3_0}-T_{3_0}N_{2_0},\;\;B_{2_0}=T_{3_0}N_{1_0}-T_{1_0}N_{3_0}\]
\[B_{3_0}=T_{1_0}N_{2_0}-T_{2_0}N_{1_0}.\]
We will, without loss of generality and up to rigid motion, choose $t_0=\frac 12,$ $T_0=(T_{1_0},T_{2_0},T_{3_0})=(1,0,0)$, and $N_0=(0,1,0).$ 
Now that we have initial conditions for $T$, $N$, and $B$, we will use the Frenet equations to express $T_0$, $T'_0$, and $T''_0$ in terms of $T_0$, $N_0$, and $B_0$ 
\[T'_0=\frac{v_0 N_0}{t_0},\;\;T''_0=(v'_0 t_0+v_0)N_0-v_0 ^2t_0 ^2T_0+v_0 ^2 \tau t_0 B_0.\]
The set of initial conditions $t_0=\frac 12,$ $T_0=(T_{1_0},T_{2_0},T_{3_0})=(1,0,0)$, and $N_0=(0,1,0)$ yields $T'_0=(0,\frac{4}{\sqrt{3} \tau},0)$ and $T''_0=(\frac{-16}{3 \tau^2},\frac{-16}{3 \sqrt{3} \tau},\frac{8}{3 \tau}).$

\subsection{Solving for $T$ via Hypergeometric functions}
$${}_p F_q(a_1,a_2, . . . , a_p;b_1,b_2, . . .,b_q;t^{a}):=\Sigma_{n=1}^{\infty}\frac{(a_1)_n . . . (a_p)_n}{(b_1)_n . . . (b_q)_n} \frac{t^{an}}{n!}$$
is the Barnes generalized hypergeometric function \cite{E}.  Note the use of the Pochhammer symbols $(x)_n:=\frac{\Gamma (x+n)}{\Gamma (x)}$.  We will also use
$${}_2 F_1^{reg}(a,b,c,t^{a}):=\frac{{}_2 F_1(a,b;c;t^{a})}{\Gamma[c]}.$$
By direct substitution (see for example section 46 of \cite{R}) it is straightforward to check that the following is a solution to (\ref{ODET}).
\[T=(T_1,T_2,T_3), \;T_j=\sum_{\ell=1}^3 c_{j\ell}S_\ell.\]
Where 
\[S_1=i t \;\; {}_3F_2(\frac 12,\frac 12,\frac 32;\frac 32-\frac i{2\tau},\frac 32+\frac i{2\tau};t^2),\]
\[S_2=(-1)^{\frac {-i}{2\tau}}t^{\frac{-i}{2\tau}} \;\; {}_3F_2(1-\frac i{2\tau},-\frac i{2\tau}, -\frac i{2\tau};\frac 12-\frac i{2\tau},1-\frac i {\tau};t^2),\]
\[S_3=(-1)^{\frac {i}{2\tau}}t^{\frac{i}{2\tau}} \;\; {}_3F_2(1+\frac i{2\tau},\frac i{2\tau}, \frac i{2\tau};\frac 12+\frac i{2\tau},1+\frac i {\tau};t^2),\]
the $c_{j\ell}$ are constants, and $i=\sqrt{-1}$.  Note that $S_1$ is pure imaginary and that $S_3$ is the complex conjugate of $S_2$.  For proper complex constants $c_{j\ell}$, $T$ is a real valued vector function.  By plugging in the initial conditions of the last section we can solve for the $c_{j\ell}$.

\subsection{Solving for $\gamma$}
\begin{figure}
\begin{center}
\includegraphics[scale=.6]{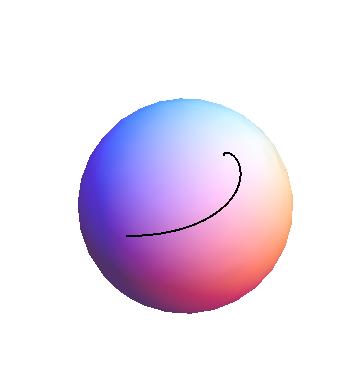}
\caption{The curve of torsion $\tau=1$ on the unit 2-sphere}  \label{curvepic}
\end{center}
\end{figure}

Recall that $\gamma(t)=\int{vT}dt.$  Since we have found $T$ in terms of hypergeometric functions, we must compute the following type of integrals:
\[\int \frac{h(t) \; {}_p F_q(a_1,a_2, . . . , a_p;b_1,b_2, . . .,b_q;,t^2)}{\tau\sqrt{1-t^2} } dt:=\int \frac{h(t) \Sigma_{n=1}^{\infty}\frac{(a_1)_n . . . (a_p)_n}{(b_1)_n . . . (b_q)_n} \frac{t^{2 n}}{n!}}{\tau\sqrt{1-t^2 }}dt\]
\[ =\frac{\alpha}{\tau}\Sigma_{n=1}^{\infty}\frac{(a_1)_n . . . (a_p)_n}{(b_1)_n . . . (b_q)_n n! } \int\frac{t^{\beta n}}{\sqrt{1-t^2 }}dt. \]
Where $\alpha, \beta$ are constants.    We repeat this process for each $S_\ell$, using the notation $\gamma=(U_1,U_2,U_3)$.
For $S_1$
\[U_1:=\int \frac{S_1}{\tau \sqrt{1-t^2}}dt
=\int \frac{i t \;\; {}_3F_2(\frac 12,\frac 12,\frac 32;\frac 32-\frac i{2\tau},\frac 32+\frac i{2\tau};t^2)}{\tau \sqrt{1-t^2}}dt\]
\[=\Sigma_{n=0}^{\infty} d_{1n}\int \frac{t^{2n+1}}{\sqrt{1-t^2}}dt ,\]
where $d_{1n}=\frac{i(1+\tau^2)\Gamma(\frac{1}{2}+n)^2\Gamma(\frac{3}{2}+n) Sech(\frac{\pi}{2 \tau})}
{2 \sqrt{\pi} \tau^3 n! \Gamma(\frac{3}{2} +n - \frac{i}{2 \tau})
\Gamma(\frac{3}{2} +n + \frac{i}{2 \tau})}.$
For $S_2$
\[U_2:=\int \frac{S_2}{\tau\sqrt{1-t^2}}dt
=\int \frac{(-1)^{\frac {-i}{2\tau}}t^{\frac{-i}{2\tau}} \;\; {}_3F_2(1-\frac i{2\tau},-\frac i{2\tau}, -\frac i{2\tau};\frac 12-\frac i{2\tau},1-\frac i {\tau};t^2)}{\tau\sqrt{1-t^2}}dt\]
\[=\Sigma_{n=0}^{\infty} d_{2n} \int \frac{t^{2n-\frac{i}{\tau}}}{\sqrt{1-t^2}}dt, \]
where $d_{2n}=\frac{e^{\frac{\pi}{2\tau}}2^{-\frac{i}{\tau}}
\Gamma(n-\frac{i}{2 \tau})^2 \Gamma(1+n-\frac{i}{2\tau})
\Gamma(\frac{-i+\tau}{2\tau})^2}{\sqrt{\pi} \tau \Gamma(1+n)
\Gamma(1+n-\frac{i}{\tau})\Gamma(-\frac{i}{2\tau})^2\Gamma(n+\frac{-i+\tau}{2\tau})}.$
For $S_3$
\[U_3:=\int \frac{S_3}{\tau\sqrt{1-t^2}}dt
=\int \frac{(-1)^{\frac {i}{2\tau}}t^{\frac{i}{2\tau}} \;\; {}_3F_2(1+\frac i{2\tau},\frac i{2\tau}, \frac i{2\tau};\frac 12+\frac i{2\tau},1+\frac i {\tau};t^2)}{\tau\sqrt{1-t^2}}dt\]
\[=\Sigma_{n=0}^{\infty} d_{3n} \int \frac{t^{2n+\frac{i}{\tau}}}{\sqrt{1-t^2}}dt, \]
where $d_{3n}=\frac{e^{-\frac{\pi}{2\tau}}2^{\frac{i}{\tau}}
\Gamma(n+\frac{i}{2 \tau})^2 \Gamma(1+n+\frac{i}{2\tau})
\Gamma(\frac{i+\tau}{2\tau})^2}{\sqrt{\pi} \tau \Gamma(1+n)
\Gamma(1+n+\frac{i}{\tau})\Gamma(\frac{i}{2\tau})^2\Gamma(n+\frac{i+\tau}{2\tau})}.$
Once again we are lucky and for each $U_\ell$ we can evaluate the integrals.  In each case they are hypergeometric.
\[U_1=\Sigma_{n=0}^{\infty} \frac{n!}{2}d_{1n}\; 
{}_2 F_1^{reg}(\frac{1}{2},1+n,2+n,t^2) ,\]
\[U_2=\Sigma_{n=0}^{\infty} \frac{\Gamma(n+\frac{-i+\tau}{2\tau})}{2}d_{2n}\; {}_2 F_1^{reg}(\frac{1}{2},n+\frac{-i+\tau}{2\tau},
\frac{3}{2}+n-\frac{i}{2\tau},t^2) ,\]
\[U_3=\Sigma_{n=0}^{\infty} \frac{\Gamma(n+\frac{i+\tau}{2\tau})}{2}d_{3n}\; {}_2 F_1^{reg}(\frac{1}{2},n+\frac{i+\tau}{2\tau},
\frac{3}{2}+n+\frac{i}{2\tau},t^2).\]
Each of the  ${}_2 F_1^{reg}$'s also has a power series.
\[{}_2 F_1^{reg}(\frac{1}{2},1+n,2+n,t^2)  = \sum_{m=0}^\infty e_{1m} t^{2m},
\;e_{1m}=\frac{\Gamma(\frac 12 +m)}{(n+m+1)\sqrt{\pi}\Gamma(1+n)\Gamma(1+m)},\]
\begin{multline*}
{}_2 F_1^{reg}(\frac{1}{2},n+\frac{-i+\tau}{2\tau},
\frac{3}{2}+n-\frac{i}{2\tau},t^2) =
 \sum_{m=0}^\infty e_{2m} t^{2m},
\\e_{2m}=\frac{2\tau\Gamma(\frac 12 +m)}
{\sqrt{\pi}\Gamma(2n\tau+2m\tau+\tau-i)\Gamma(1+m)\Gamma(n+\frac{-i+\tau}{2\tau})},
\end{multline*}
\begin{multline*} {}_2 F_1^{reg}(\frac{1}{2},n+\frac{i+\tau}{2\tau},
\frac{3}{2}+n+\frac{i}{2\tau},t^2)= \sum_{m=0}^\infty e_{3m} t^{2m},
\\e_{3m}=\frac{2\tau\Gamma(\frac 12 +m)}
{\sqrt{\pi}\Gamma(2n\tau+2m\tau+\tau+i)\Gamma(1+m)\Gamma(n+\frac{i+\tau}{2\tau})}.
\end{multline*}
Thus
\[U_1=\Sigma_{m=0}^\infty \Sigma_{n=0}^\infty \frac{n!}{2}d_{1n} e_{1m} t^{2m+2n+2},\]
\[U_2=\Sigma_{m=0}^\infty \Sigma_{n=0}^\infty \frac{\Gamma(n+\frac{-i+\tau}{2\tau})}{2} d_{2n} e_{2m} t^{2m+2n+2},\]
\[U_3=\Sigma_{m=0}^\infty \Sigma_{n=0}^\infty \frac{\Gamma(n+\frac{i+\tau}{2\tau})}{2} d_{3n} e_{3m} t^{2m+2n+2}.\]
These complicated double sums combine nicely and simplify as follows.
\[U_1=\frac{i}{2\sqrt{\pi}\tau}\Sigma_{k=0}^\infty \frac{\Gamma(\frac12 +k)}{\Gamma(2 +k)}\;{}_4F_3 (\frac12,\frac12,\frac32,-k;\frac12-k,\frac32-\frac{i}{2 \tau}, \frac32+\frac{i}{2 \tau};1)\; t^{2+2k},\]
\begin{multline*}
U_2=\frac{e^{\frac{\pi}{2 \tau}}}{\sqrt{\pi}} \Sigma_{k=0}^\infty \frac{\Gamma(\frac12 +k)}{(-i+(1+2k)\tau) \Gamma(1+k)}\;\; \times
\\ {}_4F_3 (-k,1-\frac{i}{2 \tau},-\frac{i}{2 \tau},-\frac{i}{2 \tau};\frac12-k,\frac12-\frac{i}{2 \tau},1-\frac{i}{\tau};1)
\; t^{1-\frac{i}{\tau}+2k},
\end{multline*}
\begin{multline*}
U_3=\frac{e^{-\frac{\pi}{2 \tau}}}{\sqrt{\pi}} \Sigma_{k=0}^\infty \frac{\Gamma(\frac12 +k)}{(i+(1+2k)\tau) \Gamma(1+k)}\;\; \times
\\ {}_4F_3 (-k,1+\frac{i}{2 \tau},\frac{i}{2 \tau},\frac{i}{2 \tau};\frac12-k,\frac12+\frac{i}{2 \tau},1+\frac{i}{\tau};1)
\; t^{1+\frac{i}{\tau}+2k}.
\end{multline*}
Thus we can write $\gamma=(\gamma_1,\gamma_2,\gamma_3)$ as a power series in $t$ where
\[\gamma_j=\sum_{\ell=1}^3 c_{j\ell}U_\ell.\]
The curve with $\tau=1$ is given in Figure \ref{curvepic}, this time using the explicit formula.

%



\section{Appendix}
The purpose of this Appendix is to address questions and issues about the curves raised by the referees and the editor.

For $\tau=0$, the curves are also planar and are precisely the set of circles lying on the sphere.  If we consider curves corresponding to solutions (\ref{spherekappa}) with $C=0$, and  $k(\frac{\frac{\pi}{2}}{\tau}) = \csc{\frac{\frac{\pi}{2}}{\tau} \tau}=1$,  then as $\tau$ varies from $0$ to $\infty$ the curves numerically appear to vary (in a non-uniform way) from an infinitely covered great circle, through a family of spiral ``clothoid" like curves.  

\begin{figure}
\begin{center}
\includegraphics[scale=.6]{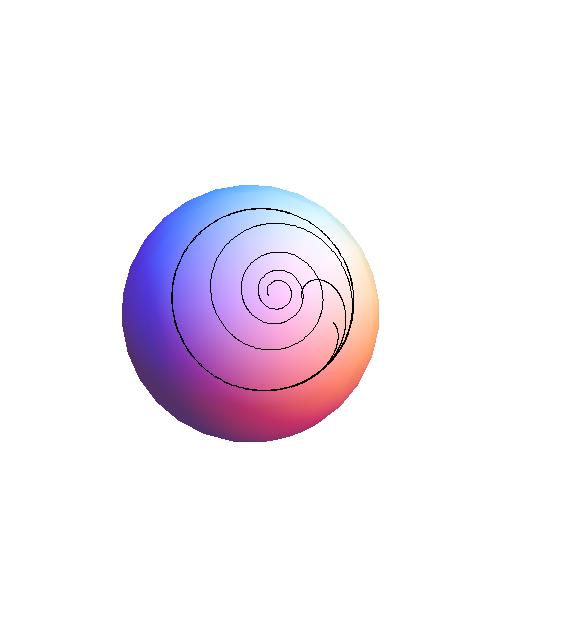}
\caption{Curves of constant torsion approaching a circle }\label{finn}
\end{center}
\end{figure}

The editor pointed out that if we consider those solutions to (\ref{spherekappa}) with $C=0$, and $k(s_0)=\csc(s_0 \tau) > 1$, then the corresponding curves approach a ``small circle" on $S^2(1)$ of constant curvature $k(s_0)$.  This is an interesting example of non-uniform convergence.  The curves as a whole converge pointwise to an infinitely covered great circle, while it is still possible to find sequences of ``tails" that converge to infinitely covered small circles.   This phenomenon is indicated in Figure \ref{finn} where one sees a sequence of curves converging to a small circle.  More details of this simple yet interesting behavior will be written up elsewhere.

We will mention a few of the qualitative properties of these curves.  Let us consider the case of curves in $S^2(1)$ with a fixed initial point and varying $\tau$.  All curves of constant torsion differ from one of these by a rigid motion. In \cite{C}, p185, it is shown that the curves are embedded, spiral infinitely often about a limiting endpoint, and are reflectionally symmetric through the initial point.  (In Figure \ref{peter} we show only the upper-half of the curves.)    \cite{C}, p185, also shows that as $\tau$ varies from $0$ to $\infty$, the length varies from $\infty$ to $0$, see also equation (\ref{spherekappa}). 

One referee asked if it would be possible to foliate $S^2(1)$ with curves of constant torsion (other than by the just using circles).  It may be possible to foliate $S^2(1)$, in some convoluted way, by packing $S^2(1)$ with pieces of curves of constant torsion; however our conjecture would be that is not possible to foliate in any ``reasonable" way.  The reasoning is as follows.  It seems to be a difficult problem to find an explicit formula for the upper endpoint in terms of the $\tau$ and the initial point.  Nevertheless, numerically as $\tau$ varies from $0$ to $\infty$ the  upper endpoint steadily moves downward from the north pole to the initial point.  In particular this would imply that the curves corresponding to an infinitesimal change in $\tau$ would (repeatedly) intersect.  It would follow that any foliation of the $S^2(1)$ by curves of constant torsion would have to include curves with common endpoints that differ by a rigid motion; a rotation about the upper endpoint.  This type of foliation could only work in some radius about the upper endpoint, because the effect of a rotation on the opposite lower endpoint would result in (repeated) intersections.  In summary, the numerics strongly indicate that there is no foliation (singular or not) of $S^2(1)$ by curves of constant torsion.

Weiner \cite{W} proved that there exist arbitrarily short closed constant torsion curves in $\mathbb{R}^3$. More recently, Musso \cite{M} studied those curves of constant torsion in $\mathbb{R}^3$, whose normal vectors sweep out elastic curves in $S^2(1)$.  Ivey \cite{I} generalizes Musso's results and gives examples of closed constant torsion curves of various knot types.  The examples in the current paper complement these known examples.

\end{document}